\documentclass[12pt]{article}
\usepackage{graphicx} 
\usepackage[english]{babel}

\usepackage[utf8]{inputenc}
\usepackage[left=1in,top=1in,right=1in,bottom=1in]{geometry}

\usepackage{amsmath}
\usepackage{graphicx}
\usepackage[colorlinks=true, allcolors=blue]{hyperref}
\usepackage{amsmath, amssymb, amsthm, amsfonts}
\usepackage{fullpage}
\usepackage{hyperref}
\usepackage{enumitem}
\usepackage{mathtools}
\usepackage[capitalise]{cleveref}

\newtheorem{theorem}{Theorem}

\newtheorem{corollary}{Corollary}

\DeclarePairedDelimiter\norm{\lVert}{\rVert}%
\DeclarePairedDelimiter\abs{\lvert}{\rvert}%
\makeatletter
\let\oldabs\abs
\def\abs{\@ifstar{\oldabs}{\oldabs*}}
\let\oldnorm\norm
\def\norm{\@ifstar{\oldnorm}{\oldnorm*}}
\makeatother

\newcommand{\ex}{{\rm ex}}

\title{Hypergraphs without complete partite subgraphs}
\author{
Dhruv Mubayi\thanks{Department of Mathematics, Statistics and Computer Science, University of Illinois, Chicago, IL 60607. Email: mubayi@uic.edu. Research partially supported by NSF Award DMS-2153576.}}
\date{\today}

\begin{document}

\maketitle
\begin{abstract}
Fix  integers $r \ge 2$ and $1\le s_1\le \cdots \le s_{r-1}\le  t$ and set $s=\prod_{i=1}^{r-1}s_i$. Let $K=K(s_1, \ldots, s_{r-1}, t)$ denote the complete $r$-partite $r$-uniform hypergraph with parts of size $s_1, \ldots, s_{r-1}, t$.
We prove that the Zarankiewicz number $z(n, K)= n^{r-1/s-o(1)}$ provided $t> 3^{s+o(s)}$.
Previously this was known only for $t > ((r-1)(s-1))!$ due to Pohoata and Zakharov. Our novel approach, which uses Behrend's construction of sets with no 3 term arithmetic progression, also applies for small values of $s_i$, for example,   it gives  $z(n, K(2,2,7))=n^{11/4-o(1)}$ where the exponent 11/4 is optimal, whereas previously this was only known with 7 replaced by 721.
\end{abstract}

\section{Introduction}

Write $K=K(s_1, \ldots, s_r)$ for the complete $r$-partite $r$-uniform hypergraph (henceforth $r$-graph) with parts of size $s_1\le s_2\le  \cdots \le s_r$. More precisely, the vertex set of $K$ comprises disjoint sets $S_1, \ldots, S_r$, where $|S_i|=s_i$ for $1\le i \le r$, and the edge set of $K$ is 
$$\{ \{x_i, \ldots, x_r\}: (x_1, \ldots, x_r) \in S_1 \times \cdots \times S_r\}.$$
Given $K$ as above, write $\ex(n, K)$ for the maximum number of edges in an $n$-vertex $r$-graph that contains no copy of $K$ as a subhypergraph. Similarly, write $z(n, K)$ for the maximum number of edges in an $r$-partite $r$-graph $H$ with parts $X_1, \ldots, X_r$, each of size $n$, such that there is no copy of $K(s_1, \ldots, s_r)$ in $H$ with $S_i \subset X_i$ for all $1\le i \le r$ (there could be copies of $K$ in $H$, where for some $i$, $S_i \not\subset X_i$). Determining $\ex(n, K)=\ex_r(n, K)$ is usually called the Tur\'an problem, while determining  $z(n, K)=z_r(n, K)$ is called the Zarankiewicz problem (we will omit the subscript $r$ if it is obvious from context). These are fundamental questions in combinatorics with applications in analysis~\cite{AL, CCW}, number theory~\cite{OS}, group theory~\cite{L}, geometry~\cite{FPS}, and computer science~\cite{BGKRSW}.

A basic result in extremal hypergraph theory, due to Erd\H os~\cite{Erdos1964}, is the upper bound
\begin{equation} \label{kstupper}\ex(n, K(s_1, \ldots, s_r)) = O(n^{r-1/s}),\end{equation}
where $s=s_1s_2\cdots s_{r-1}$ (and, as before $s_1\le s_2\le \cdots \le s_{r-1}\le s_r$). Here $s_1, \ldots, s_r$ are fixed and the asymptotic notation is taken as $n \to \infty$. As $r$ is fixed, and $z(n,K(s_1, \ldots, s_r))\le \ex(rn, K(s_1, \ldots, s_r)) $, the same upper bound as (\ref{kstupper}) holds for $z(n,K(s_1, \ldots, s_r))$.

A major problem in extremal (hyper)graph theory is to obtain corresponding lower bounds to (\ref{kstupper})  (or prove that no such lower bounds exist). In fact, it was conjectured in~\cite{Mub} that the exponent $r-1/s$ in (\ref{kstupper}) is optimal. This question has been studied for graphs since the 1930s, and results of Erd\H os-R\'enyi and Brown~\cite{Brown} gave optimal (in the exponent) lower bounds for $K(2,t)$ and $K(3,t)$. The first breakthrough for arbitrary $s_1$ occurred in the mid 1990's by Kollar-Ronyai-Szabo~\cite{KRS} and then Alon-Ronyai-Szabo~\cite{ARS}, who proved that $\ex(n, K(s_1,s_2)) = \Omega(n^{2-1/s_1})$ as long as $s_2>(s_1-1)!$. More recently, in another significant advance,
Bukh~\cite{Bukh} has proved the same lower bound as long as $s_2 > 9^{s_1+o(s_1)}$.

For $r\ge 3$,  the first nontrivial constructions that were superior to the bound given by the probabilistic deletion method were provided in the cases $s_1=\cdots =s_{r-2}=1$ and  $K(2,2,3)$ by the current author~\cite{Mub} and, soon after for $K(2,2,2)$ by Katz-Krop-Maggioni~\cite{KKM} (see also~\cite{CPZ} for recent results on the $r$-uniform case $K(2,\ldots, 2)$ that are superior to the probabilistic deletion bound but not optimal in the exponent).
Later, optimal bounds for both $\ex(n, K)$ and $z(n, K)$  were provided by Ma, Yuan, Zhang~\cite{MYZ}  (and independently by Verstra\"ete) by extending the method of Bukh, however, the threshold for $s_r$ for which the bound holds was not even explicitly calculated. More recently, 
lower bounds matching the exponent $r-1/s$ from (\ref{kstupper}) have been 
proved  for $s_r > ((r-1)(s-1))!$ by Pohoata and Zakharov~\cite{PZ}. Here we improve this lower bound on $s_r$  substantially in the Zarankiewicz case, from factorial to exponential  at the expense of a small $o(1)$ error parameter in the exponent.  The following is our main result.

\begin{theorem} \label{main}
	Fix $r \ge 3$, and positive integers $s_1, \ldots, s_{r-1}, t$. Then as $n \to \infty$,
	$$z_r(n, K(s_1, \ldots, s_{r-1}, t)) > n^{1-o(1)} \cdot z_{r-1}(n, K(s_1, \ldots, s_{r-3},  s_{r-2}s_{r-1}, t)).$$
	\end{theorem}
Applying Theorem~\ref{main} repeatedly (or doing induction on $r$) yields
	$$z_r(n, K(s_1, \ldots, s_{r-1}, t)) > n^{r-2-o(1)} \cdot z_{2}(n, K(s, t))$$
	where $s=s_1\cdots s_{r-1}$.
Bukh~\cite{Bukh} proved that $z(n, K(s, t)) = \Omega (n^{2-1/s})$ provided $t > 3^{s+o(s)}$ and this yields the following corollary. 

\begin{corollary} Fix $r \ge 2$, and  integers $1\le s_1 \le \cdots \le  s_{r-1}<t$  where $t>3^{s+o(s)}$ and $s=s_1\cdots s_{r-1}$. Then as $n \to \infty$,
	 $$z_r(n,  K(s_1, \ldots, s_{r-1}, t)) = n^{r-1/s-o(1)}.$$
	\end{corollary}

We remark that Theorem~\ref{main} can also be applied for small values of $s_i$. For example, using the result of Alon-R\'onyai-Szab\'o~\cite{ARS} that $z(n, K(4,7)) =\Omega(n^{7/4})$, it gives 
$$z(n, K(2,2,7)) > n^{1-o(1)} \,z(n, K(4, 7)) > n^{1-o(1)} \, n^{7/4} = n^{11/4-o(1)},$$ 
where the exponent $11/4$ is tight. For contrast, the previous best result due to Pohoata and Zakharov~\cite{PZ} yields only $z(n, K(2,2,721)) > \Omega(n^{11/4})$. If, as is widely believed, $z(n, K(4,4))=\Omega(n^{7/4})$, then this would imply via Theorem~\ref{main}, that $z(n,K(2,2,4))=n^{11/4-o(1)}$. 

\section{Proof}
Write $e(H)=|E(H)|$ for a hypergraph $H$.
To prove Theorem~\ref{main}, we need the following well-known consequence of Behrend's construction~\cite{B} of a subset of $[n]$ with no 3-term arithmetic progression (see, e.g., \cite{RSz}). There exists a bipartite graph $G$ with parts of size $n$ and $n^{2-o(1)}$ edges whose edge set is a union of $n$ induced matchings.  More precisely, there are pairwise disjoint matchings $M_1, \ldots, M_n$ such that $E(G)=\cup_{i=1}^n M_i$ and for all $i,j$ the edge set $M_i \cup M_j$ contains no path with three edges. Additionally, 
$e(G) =\sum_i|M_i|=  n^{2-o(1)}$.
\bigskip

{\bf Proof of Theorem~\ref{main}.} Let $H'$ be an $(r-1)$-partite $(r-1)$-graph with parts $X_1, \ldots, X_{r-3}, [n]$ and $Y $ each of size $n$ with $e(H')=z_{r-1}(n, K')$ that contains no copy of the complete $(r-1)$-partite $(r-1)$-graph $K'= K(s_1, \ldots, s_{r-3},  s_{r-2}s_{r-1}, t)$. Here we only assume that there are no copies of $K'$ where the $i$th part  of size $s_i$ is a subset of $X_i$ for $1\le i\le r-3$, the $r-2$ part of size $s_{r-2}s_{r-1}$ is a subset of $[n]$, and the $r$th part of size $t$ is a subset of $Y$. Write $d(j)$ for the degree in $H'$ of vertex $j \in [n]$, so $e(H') = \sum_j d(j)$. By relabeling we may assume that $d(1) \ge d(2) \ge \cdots \ge d(n)$.

Let $G$ be a bipartite graph with parts $A$ and $B$, each of size $n$ comprising $n$ induced matchings $M_1, \ldots, M_n$, with $e(G)=\sum_i|M_i| 
=n^{2-o(1)}$. Moreover, we may assume that $|M_1| \ge |M_2| \ge \cdots \ge |M_n|$. 

Now define the $r$-partite $r$-graph $H$ as follows: the parts of $H$ are $X_1, \ldots, X_{r-3}, A, B, Y$, each of size $n$. For each $j \in [n]$, let 
$$E_j = \{\{x_1,\ldots, x_{r-3}, a, b, y\}: \{x_1,\ldots, x_{r-3}, j, y\} \in E(H'), (a,b) \in A \times B, \{a,b\} \in M_j\}$$
and let  $E(H) = \cup _{j=1}^n E_j$. Observe that $e(H) = \sum_j d(j)|M_j|$.
 In  words, we have replaced vertex $j$ that lies in edge $\{x_1,\ldots, x_{r-3}, j, y\}$ of $H'$ by all possible pairs $ab$ of $M_j$  to create $|M_j|$  edges of $H$. Now  Chebyshev's sum inequality and $e(G)=n^{2-o(1)}$ yield
$$\frac{1}{n} \sum_{j=1}^n d(j)|M_j| \ge \frac{1}{n^2} \sum_{j=1}^n d(j) \sum_{j=1}^n |M_j| = \frac{1}{n^2} e(H') e(G)=e(H') n^{-o(1)}.$$
Hence $e(H) =\sum_j d(j)|M_j| = n^{1-o(1)}e(H')$ as required.
	
	Now suppose there is a copy $L$ of $K=K(s_1, \ldots, s_{r-1}, t)$ in $H$ where the part of size $s_i$ lies in $X_i$ for $1\le i \le r-3$, the part $A'$ of size $s_{r-2}$ lies in $A$, the $B'$ part of size $s_{r-1}$ lies in $B$, and the part of size $t$ lies in $Y$. Then all $s_{r-2}s_{r-1}$ pairs in ${V(L)\choose 2}$ within $A' \times B'$ must come from different matchings $M_i$ as the matchings are induced. Indeed, if there is an $i$ such that  $ab$ and $a'b'$ are distinct edges of $M_i$, where $a,a' \in A'$ and $b, b' \in B'$, then $ab'$ cannot lie in any edge of $H$, as $M'$ is an induced matching, but $ab'$ must lie in many edges of $L$, contradiction. The number of these matchings $M_i$ is therefore $|A'||B'|=s_{r-2}s_{r-1}$ and each such matching $M_j$ corresponds to a vertex $j$ of $[n]$.  This means that we have a forbidden copy of $K'$ in $H'$, contradiction.   
\qed
\bigskip

{\bf Remarks.}

\begin{itemize}
	
	\item One shortcoming of our  approach is that it applies only to the Zarankiewicz problem and not the Tur\'an problem. It would be interesting to rectify this.  
	
	\item For some $r$-partite $r$-graphs $H$ one can define an appropriate $(r-1)$-partite $(r-1)$-graph $H'$ such that $z_r(n, H)> n^{1-o(1)} z_{r-1}(n, H')$. This may give some further new results for hypergraphs.  
	\end{itemize}


\begin{thebibliography}{11}
\bibitem{AL} S. V. 	Astashkin, K. V.  Lykov, 
Sparse Rademacher chaos in symmetric spaces.(Russian. Russian summary)
Algebra i Analiz 28 (2016), no. 1, 3--31; translation in
St. Petersburg Math. J. 28 (2017), no. 1, 1--20

	\bibitem{ARS} N. Alon, L. R\'onyai,  T. Szab\'o, Norm-graphs: variations and applications, J. Combin. Theory Ser. B 76 (1999), no. 2, 280--290.


\bibitem{BGKRSW} L.	Babai, A. G\'al, J. Koll\'ar, L. R\'onyai,  T. Szab\'o,  A. Wigderson, 
Extremal bipartite graphs and superpolynomial lower bounds for monotone span programs, Proceedings of the Twenty-eighth Annual ACM Symposium on the Theory of Computing (Philadelphia, PA, 1996), 603--611.


\bibitem{B} Behrend, F. A, On sets of integers which contain no three terms in arithmetical progression, Proc. Nat. Acad. Sci. U.S.A. 32 (1946), 331–332.

\bibitem{Brown} W. G. Brown, On graphs that do not contain a Thomsen graph, Canad. Math. Bull. 9 (1966),
281--285.

\bibitem{Bukh} B. Bukh, 
Extremal graphs without exponentially small bicliques.
Duke Math. J. 173 (2024), no. 11, 2039--2062.

\bibitem{CCW}	A. Carbery, M. Christ, J.  Wright, 
Multidimensional van der Corput and sublevel set estimates.
J. Amer. Math. Soc. 12 (1999), no. 4, 981--1015.


\bibitem{CPZ} D. Conlon, C.  Pohoata, D. Zakharov, 
Random multilinear maps and the Erdős box problem.
Discrete Anal. 2021, Paper No. 17, 8 pp.

\bibitem{Erdos1964}	P. Erd\H os, 
On extremal problems of graphs and generalized graphs.
Israel J. Math. 2 (1964), 183--190.




\bibitem{FPS} J. Fox, Jacob, J. Pach, A. Suk, 
A polynomial regularity lemma for semialgebraic hypergraphs and its applications in geometry and property testing.
SIAM J. Comput. 45 (2016), no. 6, 2199--2223.


\bibitem{KRS}	J. Koll\'ar, L. R\'onyai, and T. Szab\'o, Norm-graphs and bipartite Turán numbers, Combinatorica 16 (1996), no. 3, 399--406.

\bibitem{KKM} N. H. Katz, E. Krop, M. Maggioni, 
Remarks on the box problem.
Math. Res. Lett. 9 (2002), no. 4, 515--519.

\bibitem{L}
A. 	Lucchini, 
Applying extremal graph theory to a question on finite groups.
Acta Math. Hungar. 168 (2022), no. 1, 295--299.

\bibitem{MYZ} J. Ma, X.  Yuan, M. Zhang, 
Some extremal results on complete degenerate hypergraphs.
J. Combin. Theory Ser. A 154 (2018), 598--609.

\bibitem{Mub} D. Mubayi,
Some exact results and new asymptotics for hypergraph Turán numbers.
Combin. Probab. Comput. 11 (2002), no. 3, 299--309.

\bibitem{OS} A.	Oppenheim, M.  Shusterman, 
Squarefree polynomials with prescribed coefficients.
J. Number Theory 187 (2018), 189--197.


\bibitem{PZ} 	C. Pohoata, D. Zakharov, Norm hypergraphs, Combinatorica, to appear. 


\bibitem{RSz}  Ruzsa, I. Z.  and Szemer\'edi, E,  Triple systems with no six points carrying three triangles, in Combinatorics (Proc. Fifth Hungarian Colloq., Keszthely, 1976), Vol. II, Colloq. Math. Soc. J\'anos Bolyai, vol. 18, North-Holland, Amsterdam-New York, 1978, 939--945.


\end{thebibliography}
\end{document}